\documentclass[12pt]{article}

\usepackage{    amsmath,
                amsfonts,
                amssymb,
                amsthm,
		graphics,
                pvDstyle,       %% my style
                verbatim,
}

\InputIfFileExists{diagrams.tex}{}{%
\IfFileExists{diagrams.sty}{\usepackage{diagrams}}
{\message{Without diagrams you can not process this file!}}}
\diagramstyle[height=2em,balance,righteqno,PostScript=dvips]

\newcommand{\rotate}[1]{\rotatebox{90}{#1}}

\newcommand\LCS{{\mathcal{LCS}h}}
\newcommand\Sh{{\mathcal{S}h}}

\newarrow{Mono}{boldhook}---{->}
\newarrow{TTo}----{->}

\newcommand\CC{{\mathbb C}}
\newcommand\ZZ{{\mathbb Z}}

\newcommand{\ca}{{\mathcal A}}
\newcommand{\cb}{{\mathcal B}}

\newcommand{\cf}{{\mathcal F}}
\newcommand{\cs}{{\mathcal S}}

\newcommand{\cu}{{\mathcal U}}

\newcommand{\cw}{{\mathcal W}}
\newcommand{\cx}{{\mathcal X}}
\newcommand{\cy}{{\mathcal Y}}
\newcommand{\cz}{{\mathcal Z}}

\newcommand{\D}{{\mathfrak D}}

\newcommand{\Dbc}{D^{b,c}}
\newcommand{\0}{\phantom0}
\newcommand{\bull}{{\scriptscriptstyle\bullet}}
\newcommand{\shortto}{\mskip-1mu\to\mskip-1mu}

\newcommand{\lococit}{loc.\ cit.}
\newcommand{\coco}{\nobreakdash-\hspace{0pt}cohomologically constructible}
\newcommand{\clc}{\nobreakdash-\hspace{0pt}cohomologically locally constant}
\newcommand{\n}[1]{\nobreakdash-\hspace{0pt}}

\let\<\langle
\let\>\rangle

\let\dlgn\boxtimes
\newcommand{\Ddlgn}{{\boxtimes^D}}
\let\ge\geqslant
\let\le\leqslant
\let\tens\otimes

\DeclareMathOperator\Ext{Ext}
\newcommand{\ExtY}{\vphantom{E}^Y{\operatorname{Ext}}}
\DeclareMathOperator\Hom{Hom}
\newcommand{\HOM}{\underline{\operatorname{Hom}}}
\DeclareMathOperator\id{id}
\DeclareMathOperator\Id{Id}
\DeclareMathOperator\modul{-mod}
\DeclareMathOperator\Ob{Ob}
\newcommand{\op}{{\operatorname{op}}}
\DeclareMathOperator\Perv{Perv}
\DeclareMathOperator\pr{pr}
\newcommand{\RHom}{R\operatorname{Hom}^\bull}
\newcommand{\RHOM}{R\underline{\operatorname{Hom}}}
\DeclareMathOperator\SimpPerv{SimpPerv}

\DeclareMathOperator\Vect{-Vect}

\newcommand{\corref}[1]{Corollary~\ref{#1}}
\newcommand{\defref}[1]{Definition~\ref{#1}}
\newcommand{\figref}[1]{Fig.~\ref{#1}}
\newcommand{\lemref}[1]{Lemma~\ref{#1}}
\newcommand{\propref}[1]{Proposition~\ref{#1}}

\newcommand{\secref}[1]{Section~\ref{#1}}

\usepackage[mathcal]{euscript}  %% redefines mathcal

\begin{document}
\bibliographystyle{amsplain}

\title{External tensor product of categories of perverse sheaves}
\author{Volodymyr Lyubashenko}

 \begin{comment}
\institute{
Institute of Mathematics,
Ukrainian National Academy of Sciences,
3, Tereshchenkivska st.,
Kyiv-4, 252601,
Ukraine;
\email{lub@imath.kiev.ua}}
\subclass{Primary 18E30}       %% ???????????
 \end{comment}

\date{November 26, 1999; math.CT/9911207}
% \date{Received: October 6, 1999 / Revised version: \today}

\maketitle

%\begin{comment}
\begin{center}
Institute of Mathematics \\
Ukrainian National Academy of Sciences \\
3, Tereshchenkivska st. \\
Kyiv-4, 252601 \\
Ukraine \\
{\catcode`@=11
lub@imath.kiev.ua}
\end{center}
%\end{comment}

\begin{abstract}
Under some assumptions we prove that the Deligne tensor
product of categories of constructible perverse sheaves on
pseudomanifolds $X$ and $Y$ is the category of
constructible perverse sheaves on $X\times Y$.
The Deligne external tensor product functor is identified
with the geometrical external tensor product.
\end{abstract}

\allowdisplaybreaks[1]

\section{Introduction}
% \footnotetext[1]{Research was supported in part by National Science
% Foundation grant 530666.}
The aim of this article is to show that the
geometrical external product of perverse sheaves is a
concrete realisation of their abstract external tensor product.
In more detail, there is a functor $\dlgn$, which assigns
to a pair of perverse sheaves $F$ on a space $X$ and
$G$ on a space $Y$ their geometrical external
tensor product $F\dlgn G$, which is a perverse sheaf
on the product of spaces $X\times Y$.
We claim that under certain assumptions the functor
$\dlgn$ makes the category of perverse sheaves on
$X\times Y$ into Deligne's tensor product of abelian
categories of perverse sheaves on $X$ and on $Y$.
This statement gives a little bit of support to the attempts
of finding triangulated Hopf category analogs for
quantum groups \cite{Lyu-hoca-SL2,Lyu-hoca-n+}.

Let us describe the objects we deal with. Unless otherwise
specified, all topological spaces are assumed
locally compact, locally completely paracompact,
locally contractible and of finite cohomological
dimension over $\CC$. A stratified space
for us will mean a compactifiable topological stratified
pseudomanifold $X$ with a stratification $\cx$ -- a locally
closed partition of $X$. In particular, $X$ can be compact or
a complex algebraic variety with algebraic strata. By
definition of a pseudomanifold there exists a filtration by
closed subspaces
\[ \cf_X: \quad X=X_n\supset X_{n-1} \supset\dots\supset
X_1\supset X_0=\varnothing, \]
such that the $S_i=X_i-X_{i-1}$ are topological manifolds,
topologically disjoint union of strata and the conditions of
\cite[Definition I.1.1]{Borel-Int_Cohom-I} hold.

Sheaves will be the sheaves of $\CC$\n-vector spaces.
Following Borel \cite[V.3.3]{Borel-Int_Cohom-V} we say that a
complex of sheaves $K$ is $\cx$\clc\ if $H^\bull(K)$
are locally constant on each stratum.
We say that $K$ is $\cx$\coco\ if it is $\cx$\clc\ and
the stalks of $H^\bull(K)$ are finite dimensional.
The full subcategory of the bounded derived
category $D^b(X)$ consisting of $\cx$\clc\ complexes is
denoted $D^b_\cx(X)$, its subcategory consisting of
$\cx$\coco\ complexes is denoted $\Dbc_\cx(X)$.

Assume that $(X,\cx)$ is equipped with a function
$p:\cx\to\ZZ$, the perversity, satisfying the condition
\[ p(S) \ge p(T) \text{ if } S \subset \overline{T}. \]
Beilinson, Bernstein and Deligne associate a
$t$\n-structure on $\Dbc_\cx(X)$ with the perversity:
the full subcategory $\vphantom{D}^pD^{\le0}_\cx(X)$
(resp. $\vphantom{D}^pD^{\ge0}_\cx(X)$) formed by
complexes $K\in\Dbc_\cx(X)$ such that for each
stratum $i_S:S \rMono X$ the following holds:
$H^mi_S^*K=0$ for $m>p(S)$ (resp. $H^mi_S^!K=0$ for
$m<p(S)$). The complexes that satisfy the both conditions
 are called perverse sheaves. The category of perverse
sheaves, the heart, is denoted
\[ \Perv(X) = \Perv(X,\cx,p) = \vphantom{D}^pD^{0}_\cx(X) =
\vphantom{D}^pD^{\le0}_\cx(X) \cap \vphantom{D}^pD^{\ge0}_\cx(X). \]

The external tensor product functor
\[ \dlgn: \Dbc_\cx(X)\times\Dbc_\cy(Y) \to \Dbc_{\cx\times\cy}(X\times Y) \]
is defined on $K\in\Dbc_\cx(X)$, $M\in\Dbc_\cy(Y)$ as
\[ K\dlgn M = (\pr_X^*K)\overset{L}\tens_\CC(\pr_Y^*M)
= (\pr_X^*K)\tens_\CC(\pr_Y^*M), \]
where $\pr_X: X\times Y\to X$, $\pr_Y: X\times Y\to Y$
are the projections. We will abuse the notation
denoting the derived functors in the same way as for sheaves,
the prefixes $R$ and $L$ will be often omitted.
Our main result is the following.

\begin{subtheorem}
The restriction of the external tensor product
functor to perverse sheaves gives a functor
\[ \dlgn: \Perv(X,\cx,p)\times\Perv(Y,\cy,q) \to
\Perv(X\times Y,\cx\times\cy,p\dotplus q), \]
where the perversity $p\dotplus q$ is given by
$(p\dotplus q)(S\times T) = p(S)+q(T)$ for $S\in\cx$,
$T\in\cy$. This functor makes the target category
into the Deligne tensor product of abelian
$\CC$\n-linear categories $\Perv(X,\cx,p)$ and $\Perv(Y,\cy,q)$.
\end{subtheorem}

Recall that the tensor product of abelian categories
is introduced by Deligne in \cite{Del:cat}.

The rest of the paper will be devoted to the proof of this theorem.
In \secref{Presults} we prove some preliminary and
technical results, the main of which is the isomorphism
\[ \RHOM(A,C)\dlgn \RHOM(B,D) \rTTo^\sim \RHOM(A\dlgn B,C\dlgn D) \]
for cohomologically constructible complexes
(\corref{cor-RHom-box-isom}). Using it we show that $\dlgn$ is
$t$\n-exact, and that $\dlgn$
restricts to perverse sheaves (\corref{cor-Perv2-box-Perv}).
In \secref{sec-simple-pvsh} we study simple perverse
sheaves. The relationship between $\RHom$ and $\dlgn$ is
considered in \secref{sec-Rhom-box}. We reformulate our main
theorem in \secref{sec-thm-main} and prove it in a sequence
of lemmas. Appendix~A contains a list of useful formulas.

\begin{acknowledgement}
I am grateful to P.~Deligne, D.~N. Yetter, L.~Crane, and S.~A. Ovsienko
for fruitful and enlightening discussions.
Commutative diagrams in this paper are drawn with
the help of the package {\tt diagrams.tex} of Paul Taylor.
Research was supported in part by National Science
Foundation grant 530666.
\end{acknowledgement}

\section{Preliminary results}\label{Presults}

\begin{lemma}[cf. \cite{Borel-Int_Cohom-V} 10.23(2)]
\label{lem-Bor-PD5}
Let $Z$ be a locally compact, locally completely
paracompact topological space of finite cohomological dimension.
Let $A,B,C\in D^b(Z)$. Then there is a natural morphism of functors
\[ \nu: \RHOM(A,B)\tens C \to \RHOM(A,B\tens C). \]
\end{lemma}

\begin{proof}
Let us begin with sheaves $A,B,C$ on $Z$. For an
open $U\subset Z$ there is a mapping
\begin{align*}
\Hom_{Sh(U)}(A\big|_U,B\big|_U)\times C(U) \to
\Hom_{Sh(U)}&(A\big|_U,(B\tens  C)\big|_U), \\
(f,h) \longmapsto [A(V) \to & B(V)\tens C(V) \to
(B\tens C)(V)], \\
a \mapsto & f(V)(a)\tens h\big|_V
\end{align*}
where $V\subset U$ is an arbitrary open subset.
The above mapping factorises as follows
\begin{multline*}
\Hom_{Sh(U)}(A\big|_U,B\big|_U)\times C(U) \to
\HOM(A,B)(U)\tens_\CC C(U) \to \\
\to [\HOM(A,B)\tens C](U) \rTTo^{\exists\nu(U)} \HOM(A,B\tens C)(U) \\
= \Hom_{Sh(U)}(A\big|_U,(B\tens  C)\big|_U)
\end{multline*}
by the universal property of tensor products.
So we get a sheaf morphism
\[ \nu: \HOM(A,B)\tens C \to \HOM(A,B\tens C). \]
We extend it to complexes of sheaves $A,B,C\in D^b(Z)$
without additional signs since we work with conventions
in which $A\tens\HOM^\bull(A,B)\to B$ is a chain map.
Assuming $B$ to be a complex of injective
sheaves we get the required
\begin{multline*}
\nu: \RHOM(A,B)\tens C = \HOM(A,B)\tens C \rTTo^\nu \\
\to \HOM(A,B\tens C) \to \RHOM(A,B\tens C).
\end{multline*}
\end{proof}

\begin{corollary}
For $A,B,C,D\in D^b(Z)$ there is an iterated morphism
\begin{multline*}
\RHOM(A,C)\tens \RHOM(B,D) \rTTo^\nu
\RHOM(A,C\tens \RHOM(B,D)) \\
\rTTo^{\RHOM(A,\nu')} \RHOM(A,\RHOM(B,C\tens D))
\rTTo^{\eqref{ApP10.2}} \RHOM(A\tens B,C\tens D),
\end{multline*}
where
\begin{multline*}
\nu': C\tens\RHOM(B,D) \rTTo^\sigma_\sim
\RHOM(B,D)\tens C \\
\rTTo^\nu \RHOM(B,D\tens C)
\rTTo^{\RHOM(B,\sigma)}_\sim \RHOM(B,C\tens D),
\end{multline*}
and $\sigma$ is the symmetry in the category of
complexes -- the signed permutation.
\end{corollary}

\begin{corollary}\label{cor-RHom-box-exists}
For $Z=X\times Y$, $A,C\in D^b(X)$, $B,D\in D^b(Y)$
there is a morphism
\begin{multline*}
\RHOM(A,C)\dlgn \RHOM(B,D) =
p_X^*\RHOM(A,C)\tens p_Y^*\RHOM(B,D) \\
\rTTo^{\eqref{ApP10.21}}
\RHOM(p_X^*A,p_X^*C)\tens \RHOM(p_Y^*B,p_Y^*D) \\
\to \RHOM(p_X^*A\tens p_Y^*B, p_X^*C\tens p_Y^*D) =
\RHOM(A\dlgn B,C\dlgn D),
\end{multline*}
where the morphism from the previous corollary is used.
\end{corollary}

\begin{lemma}\label{lem-f*g*}
Let $f:U\to X$, $g:V\to Y$ be continuous maps. For
$A\in D^b(X)$, $B\in D^b(Y)$ we have a natural isomorphism
\[ (f\times g)^*(A\dlgn B) \simeq f^*A\dlgn g^*B. \]
\end{lemma}

\begin{proof}
Using \eqref{ApP10.1} we get
\[ (f\times g)^*(p_1^*A\tens p_2^*B) \simeq
(f\times g)^*p_1^*A \tens (f\times g)^*p_2^*B \simeq
p_1^*f^*A \tens p_2^*g^*B. \]
\end{proof}

\begin{lemma}\label{lem-i!p*=p*i!}
Let us consider the following stratified maps
of stratified pseudomanifolds.
\begin{diagram}
S\times Y & \rTTo^{i\times1} & U\times Y \\
\dTTo<{p_S} && \dTTo>{p_U} \\
S & \rTTo^i & U
\end{diagram}
For any $\cu$\coco\ $A\in\Dbc_\cu(U)$ we have
\[ (i\times1)^!p_U^*A \simeq p_S^*i^!A. \]
\end{lemma}

\begin{proof}
Let us prove the Verdier dual of this isomorphism:
\[ (i\times1)^*p_U^!B \simeq p_S^!i^*B \]
for $B=\D A\in\Dbc_\cu(U)$.
Isomorphism \eqref{eq-Ap-p!B} gives
\[ p_U^!B \simeq B\dlgn\D_Y \qquad\text{and}\qquad
p_S^!(i^*B) \simeq (i^*B)\dlgn\D_Y. \]
Hence,
\[ (i\times1)^*p_U^!B \simeq (i\times1)^*(B\dlgn\D_Y) \simeq
(i^*B)\dlgn\D_Y \simeq p_S^!i^*B. \]
\end{proof}

\begin{proposition}\label{pro-RHom-main}
Let $(X,\cx)$ be a stratified topological pseudomanifold, and
let $Y$ be a locally compact, locally completely
paracompact, locally contractible and of finite
cohomological dimension over $\CC$. Assume that
$A\in\Dbc_\cx(X)$ is $\cx$\coco, $B\in D^b_\cx(X)$ is
$\cx$\clc, and let $C\in D^b(Y)$. Then the morphism
\[ \nu: \RHOM(p_X^*A,p_X^*B)\tens p_Y^*C \to
\RHOM(p_X^*A,p_X^*B\tens p_Y^*C). \]
from \lemref{lem-Bor-PD5} is an isomorphism.
\end{proposition}

\begin{proof}
Let
\[ \cf_X: \quad X=X_n \supset X_{n-1} \supset\dots\supset
X_1 \supset X_0 = \varnothing \]
be a closed filtration of $X$ such that the connected
components of $S_i=X_i-X_{i-1}$ are strata of $X$. Let
\[ U_i = X-X_{n-i} \]
be the complement open filtration, then
$S_{n-k}=U_{k+1}-U_k=X_{n-k}-X_{n-k-1}$. Denote certain
inclusions and projections as in the following diagram.
\begin{diagram}
X\times Y & \lMono^{J_k} & U_k\times Y & \rMono^J &
U_{k+1}\times Y & \lMono^I & S_{n-k}\times Y \\
\dTTo<{p_1} && \dTTo>{p_1} && \dTTo>{p_1} && \dTTo>{p_1} \\
X & \lMono^{j_k} & U_k & \rMono^j & U_{k+1} & \lMono^i & S_{n-k}
\end{diagram}

For a complex $K\in D^b(X\times Y)$ denote
$(K)_k=J_k^*K\in D^b(U_k\times Y)$. For $N>n$ we have
$(K)_N=K$. We want to prove by induction on $k$ that
\begin{equation}\label{eq-nu-k}
\nu_k: \RHOM((p_1^*A)_k, (p_1^*B)_k) \tens (p_2^*C)_k \to
\RHOM((p_1^*A)_k, (p_1^*B)_k \tens (p_2^*C)_k)
\end{equation}
is an isomorphism. The $k=1$ (or $n=1$) case reduces to the
simplest situation, where $X$ is a disjoint union of open strata.
Since $A$ is cohomologically locally constant with finite
dimensional cohomology stalks, it may be replaced
\emph{locally} with a complex of constant sheaves with
finite dimensional stalks. So $\HOM(p_1^*A,D)$ is isomorphic
to a sum of shifted copies of $D$ (with modified differential),
hence, it coincides with $\RHOM(p_1^*A,D)$.
Clearly, $\nu_1$ is an isomorphism.

For a complex $M\in D^b(X)$ denote
$M_k=j_k^*M\in D^b(U_k)$. Then there is an isomorphism
for $M\in D^b(X)$
\[ (p_1^*M)_k = J_k^*p_1^*M \simeq p_1^*j_k^*M
= p_1^*(M_k) = p_1^*M_k. \]

Assuming that $\nu_k$ is an isomorphism, let us prove that
$\nu_{k+1}$ is an isomorphism as well. Apply to the standard triangle
\[ I_!I^!(p_1^*B)_{k+1} \to (p_1^*B)_{k+1} \to J_*(p_1^*B)_k \to \]
the both functors in \eqref{eq-nu-k}. It gives two triangles
and a morphism between them, written down in diagram
in \figref{dia-2triangles}.

\begin{figure}
\setbox0=\hbox{%
\hspace*{0mm}
\begin{diagram}
\0 \\ \0 \\ \0 \\ \0 \\ \0 \\ \0 \\ \0 \\
\RHOM( p_1^*A_{k+1}, I_!I^!p_1^*B_{k+1}) \tens (p_2^*C)_{k+1}
& \rTTo &
\RHOM( p_1^*A_{k+1}, p_1^*B_{k+1}) \tens (p_2^*C)_{k+1}
& \rTTo &
\RHOM( p_1^*A_{k+1}, J_*p_1^*B_k) \tens (p_2^*C)_{k+1} \to \\
\dTTo<{\nu'} && \dTTo<{\nu_{k+1}} && \dTTo<{\nu''} \\
\RHOM( p_1^*A_{k+1}, I_!I^!p_1^*B_{k+1} \tens (p_2^*C)_{k+1})
& \rTTo &
\RHOM( p_1^*A_{k+1}, p_1^*B_{k+1} \tens (p_2^*C)_{k+1})
& \rTTo &
\RHOM( p_1^*A_{k+1}, J_*p_1^*B_k \tens (p_2^*C)_{k+1}) \to
\end{diagram}
\hspace*{0mm}
}
\rotate{\box0}
\caption{A morphism of two distinguished triangles}
\label{dia-2triangles}
\end{figure}

Let us prove that $\nu''$ is an isomorphism. Indeed, this
morphism is a composition of several isomorphisms:
\begin{align*}
& \RHOM( p_1^*A_{k+1}, J_*p_1^*B_k) \tens (p_2^*C)_{k+1}
\rTTo^{\eqref{ApP10.3(1)}}_\sim \\
& J_*\RHOM( p_1^*A_k, p_1^*B_k) \tens (p_2^*C)_{k+1}
\rTTo^{\eqref{ApP10.21}}_\sim \\
& J_*p_1^*\RHOM( A_k, B_k) \tens (p_2^*C)_{k+1} \rTTo_\sim \\
& J_*p_1^*\RHOM( A, B)_k \tens (p_2^*C)_{k+1} \rTTo_\sim \\
\intertext{(since $\RHOM(A,B)$ is $\cx$\clc\ by
\cite{Borel-Int_Cohom-V} Theorem 8.6, we can apply Lemma 10.22 \lococit)}
& J_*[p_1^*\RHOM( A,B)_k \tens (p_2^*C)_k]
\rTTo_\sim^{J_*\nu_k} \\
& J_*[ \RHOM( p_1^*A_k, p_1^*B_k \tens (p_2^*C)_k)]
\rTTo_\sim^{\eqref{ApP10.3(1)}} \\
& \RHOM( p_1^*A_{k+1}, J_*[p_1^*B_k \tens (p_2^*C)_k] )
\rTTo_\sim \\
\intertext{(again by Lemma 10.22 \cite{Borel-Int_Cohom-V})}
& \RHOM( p_1^*A_{k+1}, J_*p_1^*B_k \tens (p_2^*C)_{k+1}).
\end{align*}

Let us prove that $\nu'$ is an isomorphism. Indeed, this
morphism is a composition of several isomorphisms:

\begin{align*}
& \RHOM( p_1^*A_{k+1}, I_!I^!p_1^*B_{k+1}) \tens (p_2^*C)_{k+1}
\rTTo^{\eqref{ApP10.3(1)}}_\sim \\
& I_!\RHOM( I^*p_1^*A_{k+1}, I^!p_1^*B_{k+1}) \tens (p_2^*C)_{k+1}
\rTTo^{\eqref{ApP10.8(2)}}_\sim \\
& I_![\RHOM( I^*p_1^*A_{k+1}, I^!p_1^*B_{k+1}) \tens I^*(p_2^*C)_{k+1}]
\rTTo_\sim \\
\intertext{(by \lemref{lem-i!p*=p*i!} with $S=S_{n-k}$,
$U=U_{k+1}$, $A=B_{k+1}$)}
& I_![\RHOM( p_1^*i^*A_{k+1}, p_1^*i^!(B_{k+1})) \tens p_2^*C]
\rTTo_\sim^{I_!\nu} \\
\intertext{(applying $n=1$ case to $S_{n-k}$ with the trivial
filtration in place of $X$)}
& I_![\RHOM( p_1^*i^*A_{k+1}, p_1^*i^!(B_{k+1}) \tens p_2^*C)]
\rTTo_\sim^{\text{\lemref{lem-i!p*=p*i!}}} \\
& I_![\RHOM( I^*p_1^*A_{k+1}, I^!p_1^*B_{k+1} \tens I^*(p_2^*C)_{k+1})]
\rTTo_\sim^{\eqref{ApP10.3(1)}} \\
& \RHOM( p_1^*A_{k+1}, I_![I^!p_1^*B_{k+1} \tens I^*(p_2^*C)_{k+1}])
\rTTo_\sim^{\eqref{ApP10.8(2)}} \\
& \RHOM( p_1^*A_{k+1}, I_!I^!p_1^*B_{k+1} \tens (p_2^*C)_{k+1}).
\end{align*}
Since $\nu'$ and $\nu''$ are isomorphisms, so is $\nu_{k+1}$.
\end{proof}

The above proposition determines when the morphisms
$\nu$ used in \corref{cor-RHom-box-exists} are
isomorphisms. So we get

\begin{corollary}\label{cor-RHom-box-isom}
Let $A\in\Dbc_\cx(X)$ be $\cx$\coco, let $B\in\Dbc_\cy(Y)$
be $\cy$\coco, let $C\in D^b_\cx(X)$ be $\cx$\clc, and
let $D\in D^b_\cy(Y)$ be $\cy$\clc. Then the morphism
\[ \RHOM(A,C)\dlgn \RHOM(B,D) \to \RHOM(A\dlgn B,C\dlgn D) \]
is an isomorphism.
\end{corollary}

\begin{corollary}\label{cor-DADB}
Let $A\in\Dbc_\cx(X)$ be $\cx$\coco, and let
$B\in\Dbc_\cy(Y)$ be $\cy$\coco. Then
\[ \D A\dlgn\D B \simeq \D(A\dlgn B). \]
\end{corollary}

\begin{proof}
Indeed,
\begin{align*}
\D A\dlgn\D B &= \RHOM(A,\D_X) \dlgn \RHOM(B,\D_Y) \\
\rTTo_\sim^{\text{\corref{cor-RHom-box-isom}}} &
\RHOM(A\dlgn B, \D_X\dlgn\D_Y) \\
\rTTo_\sim^{\eqref{ApC10.26}} &
\RHOM(A\dlgn B, \D_{X\times Y}) \\
&= \D(A\dlgn B).
\end{align*}
\end{proof}

In addition to Lemma \ref{lem-f*g*} we have

\begin{proposition}\label{prop-!*box}
\begin{enumerate}
\renewcommand\labelenumi{(\roman{enumi})}
\item Let $A\in\Dbc_\cx(X)$, $B\in\Dbc_\cy(Y)$, and let
    $f:U\to X$, $g:V\to Y$ be stratified maps. Then
    \[ (f\times g)^!(A\dlgn B) \simeq f^!A\dlgn g^!B. \]

\item
\begin{enumerate}
\renewcommand\labelenumii{(\alph{enumii})}
    \item Let $A\in D^b(X)$, $B\in D^b(Y)$, and let
    $f:X\to U$, $g:Y\to V$ be continuous maps. Then
    \[ (f\times g)_!(A\dlgn B) \simeq f_!A\dlgn g_!B. \]

\item Furthermore, if $A\in\Dbc_\cx(X)$, $B\in\Dbc_\cy(Y)$,
and the stratified maps $f$, $g$ are proper, or complex
algebraic, or every fibre of $f$, $g$ is compactifiable, then
    \[ (f\times g)_*(A\dlgn B) \simeq f_*A\dlgn g_*B. \]
\end{enumerate}
\end{enumerate}
\end{proposition}

\begin{proof}
(i) Deduce this isomorphism applied to the objects
$\D A$, $\D B$ via \lemref{lem-f*g*}:
\begin{align*}
(f\times g)^!(\D A\dlgn\D B)
\rTTo^{\text{\corref{cor-DADB}}}_\sim &
(f\times g)^!\D(A\dlgn B) \\
\rTTo^{\eqref{ApT10.17(1)}}_\sim & \D(f\times g)^*(A\dlgn B) \\
\rTTo^{\text{\lemref{lem-f*g*}}}_\sim & \D(f^*A\dlgn g^*B) \\
\rTTo^{\text{\corref{cor-DADB}}}_\sim &
\D f^*A\dlgn \D g^*B \\
\rTTo^{\eqref{ApT10.17(1)}}_\sim & f^!\D A\dlgn g^!\D B.
\end{align*}

(ii)(a) It suffices to consider the case $g=\id$. Combined
with the similar case $f=\id$, it implies the general case.
Based on the diagram
\begin{diagram}
X & \lTTo^{p_X} & X\SWpbk\times Y & \rTTo^{p_2} & Y \\
\dTTo<f && \dTTo>{f\times1} && \dEq \\
U & \lTTo^{p_U} & U\times Y & \rTTo^{p_2} & Y
\end{diagram}
the required isomorphism is composed of the following isomorphisms
\begin{align*}
(f\times1)_!(A\dlgn B) = & (f\times1)_!(p_X^*A\tens p_2^*B) \\
\simeq & (f\times1)_!(p_X^*A\tens (f\times1)^*p_2^*B) \\
\rTTo^{\eqref{ApP10.8(2)}}_\sim & (f\times1)_!p_X^*A\tens p_2^*B \\
\rTTo^{\text{base change}}_\sim & p_U^*f_!A\tens p_2^*B \\
= & f_!A\dlgn B.
\end{align*}

(ii)(b) Is deduced from (ii)(a) using \eqref{ApT10.17(2)}
similarly to (i).
\end{proof}

\begin{proposition}
The functor
\[ \dlgn: \vphantom{D}^p\Dbc_\cx(X) \times
\vphantom{D}^q\Dbc_\cy(Y) \to
\vphantom{D}^{p\dotplus q}\Dbc_{\cx\times\cy}(X\times Y) \]
is $t$-exact.
\end{proposition}

\begin{proof}
For $K\in D^b(X)$, $L\in D^b(Y)$ we have the K\"unneth formula
\begin{align*}
H^n(K\dlgn L) &= H^n(p_1^*K\tens p_2^*L) \simeq
\oplus_{k+l=n} H^kp_1^*K\tens H^lp_2^*L \\
&\simeq \oplus_{k+l=n} p_1^*H^kK\tens p_2^*H^lL =
\oplus_{k+l=n} H^kK\dlgn H^lL.
\end{align*}
Let $A\in D_\cx^{\le p}(X)$, $B\in D_\cy^{\le q}(Y)$, and let
$S\in\cx$, $T\in\cy$ be strata. If $n>p(S)+q(T)$, then
\[ H^n(i_S\times i_T)^*(A\dlgn B) = H^n(i_S^*A\dlgn i_T^*B)
\simeq \oplus_{k+l=n} H^k(i_S^*A)\dlgn H^l(i_T^*B) = 0 \]
by Lemma \ref{lem-f*g*} and by the K\"unneth formula. Hence,
$A\dlgn B\in D^{\le p\dotplus q}_{\cx\times\cy}(X\times Y)$.

Suppose now that $A\in D_\cx^{\ge p}(X)$,
$B\in D_\cy^{\ge q}(Y)$. If $n<p(S)+q(T)$, then
\[ H^n(i_S\times i_T)^!(A\dlgn B) = H^n(i_S^!A\dlgn i_T^!B)
\simeq \oplus_{k+l=n} H^k(i_S^!A)\dlgn H^l(i_T^!B) = 0 \]
by \propref{prop-!*box} and the K\"unneth formula. Hence,
$A\dlgn B\in D^{\ge p\dotplus q}_{\cx\times\cy}(X\times Y)$.
\end{proof}

\begin{corollary}\label{cor-Perv2-box-Perv}
The restriction of $\dlgn$ to perverse sheaves gives a
$\CC$\n-bilinear functor
\[ \dlgn: \Perv(X,\cx,p)\times \Perv(Y,\cy,q) \to
\Perv(X\times Y,\cx\times\cy,p\dotplus q) \]
exact in each variable.
\end{corollary}

\begin{proof}
For a fixed $B\in\Perv(Y)$ the functor
$T: \vphantom{D}^p\Dbc_\cx(X) \to
\vphantom{D}^{p\dotplus q}\Dbc_{\cx\times\cy}(X\times Y)$,
$A\mapsto A\dlgn B$ is $t$\n-exact. Hence, the functor
\[ \vphantom{T}^pT =
\vphantom{H}^{p\dotplus q}H^0 \circ T\circ\epsilon:
\Perv(X) \to \Perv(Y), \quad A\mapsto A\dlgn B, \]
is exact by \cite[Proposition 1.3.17(i)]{BeilBD:Perverse}.
\end{proof}

\section{Simple perverse sheaves}\label{sec-simple-pvsh}
\subsection{The case of trivial stratification}
In the case of a trivial stratification $\cx=\{X\}$ of a
connected manifold $X$ the perversity $p$ is an integer
$p(X)$. We have
\begin{align*}
\Perv(X,\cx,p) &=
\{ K\in\Dbc_\cx(X) \mid H^nK=0 \text{ unless } n=p(X) \} \\
&= \LCS(X)[-p(X)]
\end{align*}
and the category of locally constant sheaves of finite rank
$\LCS(X)$ is equivalent to $\pi_1(X)\modul$ -- the category
of $\pi_1(X)$-modules finite dimensional over $\CC$.

\subsection{Intersection cohomology sheaves}
By \cite[Proposition 1.4.26]{BeilBD:Perverse} any simple
perverse sheaf on $(X,\cx,p)$ comes from a simple
$\pi_1(S)$-module for one of the strata $S$ via the functor
\begin{multline*}
\pi_1(S)\modul \rTTo^\sim \LCS(S) \rTTo^\sim \LCS(S)[-p(S)] \\
= \Perv(S,\{S\},p(S)) \rTTo^{j_{S!*}} \Perv(\overline S)
\rTTo^{i_{S*}} \Perv(X),
\end{multline*}
where $\overline S$ is the closure of $S$ and
$j_S:S \rMono \overline S$, $i_S:\overline S \rMono X$ are the
inclusions. The prolongation functor $j_{S!*}$ is defined in
\cite[Definition 1.4.22]{BeilBD:Perverse}
(see \defref{def-Ap-j!*}). Here the functor $i_{S*}$ is the
restriction of a $t$\n-exact functor
$i_{S*}:\Dbc_{\overline\cs}(\overline S) \to \Dbc_\cx(X)$,
$\overline\cs=\cx\cap\overline S$
\cite[Proposition 1.4.16]{BeilBD:Perverse}.

Now let us discuss the behaviour of so obtained perverse
sheaves with respect to $\dlgn$.

\begin{proposition}\label{pro-6boxes}
Let $S\in\cx$, $T\in\cy$ be strata of $(X,\cx,p)$ and
$(Y,\cy,q)$. Then there are functorial isomorphisms
\begin{diagram}
\pi_1(S)\modul \times \pi_1(T)\modul & \rTTo^{\tens_\CC} &
\pi_1(S\times T)\modul \\
\dTTo<\wr & \simeq & \dTTo<\wr \\
\LCS(S)\times\LCS(T) & \rTTo^\dlgn & \LCS(S\times T) \\
\dTTo<\wr & \simeq & \dTTo<\wr \\
\LCS(S)[-p(S)] \times\LCS(T)[-q(T)] & \rTTo^\dlgn &
\LCS(S\times T)[-p(S)-q(T)] \\
\dEq & = & \dEq \\
\Perv(S) \times \Perv(T) & \rTTo^\dlgn & \Perv(S\times T) \\
\dTTo<{j_{S!*}\times j_{T!*}} & \simeq & \dTTo>{(j_S\times j_T)_{!*}} \\
\Perv(\overline S) \times \Perv(\overline T) & \rTTo^\dlgn &
\Perv(\overline S\times \overline T) \\
\dTTo<{i_{S*}\times i_{T*}} & \simeq & \dTTo>{(i_S\times i_T)_*} \\
\Perv(X) \times \Perv(Y) & \rTTo^\dlgn & \Perv(X\times Y)
\end{diagram}
\end{proposition}

\begin{proof}
The essential part is to construct a functorial isomorphism
\[ (j_S\times j_T)_{!*}(A\dlgn B) \simeq j_{S!*}A\dlgn j_{T!*}B. \]
Due to the isomorphism
\[ (j_S\times j_T)_{!*} \simeq (j_S\times1)_{!*} (1\times j_T)_{!*}, \]
which follows from \cite[(2.1.7.1)]{BeilBD:Perverse}, see
\eqref{Ap2.1.7.1}, we only have to prove the following lemma.
\end{proof}

\begin{lemma}
Let $U$ be an open stratified subspace of a stratified
pseudomanifold $Z=(Z,\cz,p)$. Denote $j:U\subset Z$ the
inclusion. Let $A\in\Perv(U,p)$ and $B\in\Perv(W,q)$. Then
in $\Perv(Z\times W,p\dotplus q)$ we have a functorial isomorphism
\[ (j\times 1)_{!*}(A\dlgn B) \simeq j_{!*}A\dlgn B. \]
\end{lemma}

\begin{proof}
The left hand side is determined uniquely as a prolongation
$C$ of $A\dlgn B$ to $\Dbc_{\cz\times\cw}(Z\times W)$,
that is a complex equipped with an isomorphism
$(j\times1)^*C \simeq A\dlgn B$, such that for any stratum
$s:S\rMono Z$, that is not contained in $U$, and for any
stratum $t:T\rMono W$ we have $H^m((s\times t)^*C)=0$ for
$m\ge p(S)+q(T)$ and $H^m((s\times t)^!C)=0$ for
$m\le p(S)+q(T)$ \cite[Proposition 2.1.9]{BeilBD:Perverse},
see \defref{def-Ap-j!*}.

Let us verify these conditions for the right hand side.
Indeed, we have
\begin{equation}
(j\times1)^*(j_{!*}A\times B) \simeq j^*j_{!*}A\dlgn B \simeq A\dlgn B,
\label{eq-iso-j!*AB}
\end{equation}
\begin{align*}
H^m(s\times t)^*(j_{!*}A\times B) &\simeq H^m(s^*j_{!*}A\dlgn t^*B) \\
&\simeq \oplus_{k+l=m} H^k(s^*j_{!*}A)\dlgn H^l(t^*B) \\
&\simeq \oplus_{\substack{k+l=m\\k<p(S),l\le q(T)}}
H^k(s^*j_{!*}A)\dlgn H^l(t^*B).
\end{align*}
This vanishes if $m\ge p(S)+q(T)$. Similarly,
\begin{align*}
H^m(s\times t)^!(j_{!*}A\times B) &\simeq H^m(s^!j_{!*}A\dlgn t^!B) \\
&\simeq \oplus_{k+l=m} H^k(s^!j_{!*}A)\dlgn H^l(t^!B) \\
&\simeq \oplus_{\substack{k+l=m\\k>p(S),l\ge q(T)}}
H^k(s^!j_{!*}A)\dlgn H^l(t^!B).
\end{align*}
This vanishes if $m\le p(S)+q(T)$.

Therefore, there exists a unique isomorphism
\[ \alpha_{A,B}: (j\times 1)_{!*}(A\dlgn B) \rTTo^\sim j_{!*}A\dlgn B, \]
such that the diagram
\begin{diagram}
(j\times 1)^*(j\times 1)_{!*}(A\dlgn B) &
\rTTo_\sim^{(j\times1)^*\alpha_{A,B}} &
(j\times 1)^*(j_{!*}A\dlgn B), \\
\dTTo<\wr && \dTTo<\wr>{\eqref{eq-iso-j!*AB}} \\
A\dlgn B & \rEq & A\dlgn B
\end{diagram}
is commutative. Uniqueness of $\alpha_{A,B}$ implies its
functoriality in $A$, $B$.
\end{proof}

\section{Right derived $\Hom^\bull$ and the
external tensor product}\label{sec-Rhom-box}
The right derived $\Hom^\bull$ can be defined as the functor
\[ \RHom: D(X)^\op \times D(X) \to D(\CC\Vect) \]
equipped with an isomorphism
\begin{diagram}
K(\Sh(X))^\op \times K(\Sh(X)) & \rTTo^{\Hom^\bull} &
K(\CC\Vect) \\
\dTTo && \dTTo \\
D(X)^\op\times D(X) & \rTTo^{\RHom} & D(\CC\Vect)
\end{diagram}
Its value on $A, B\in D(X)$ can be computed as
$\Hom^\bull(A,I)$ for a quasiisomorphism $B\to I$ with a
complex $I$ of injective sheaves. The object
$\RHOM(A,B)\in D(X)$ also can be computed as $\HOM(A,I)$,
which is a complex of flabby sheaves. To compute the value
of $R\Gamma=Rp_*:D(X) \to D(\CC\Vect)$ for $p:X\to pt$ on
$\RHOM(A,B)$ we can use the flabby resolution $\HOM(A,I)$
of $\RHOM(A,B)$. Hence,
\[ Rp_*\RHOM(A,B) \simeq p_*\HOM(A,I) = \Hom^\bull(A,I)
\simeq \RHom(A,B), \]
or $R\Gamma\circ\RHOM\simeq\RHom$.

Let us apply \propref{prop-!*box}(ii)(b) to maps
$p_X:X\to pt$, $p_Y:Y\to pt$. Recall that the stratified
pseudomanifolds $X$, $Y$ are assumed in this paper to be
compactifiable (for instance, compact or complex algebraic).

\begin{lemma}\label{lem-RGamma}
For arbitrary $E\in\Dbc_\cx(X)$ and $F\in\Dbc_\cy(Y)$ we have
\[ R\Gamma(E\dlgn F) \simeq R\Gamma E\tens_\CC R\Gamma F. \]
\end{lemma}

\begin{proof}
$(p_X\times p_Y)_*(E\dlgn F) \simeq p_{X*}E\dlgn p_{Y*}F =
p_{X*}E\tens_\CC p_{Y*}F$.
\end{proof}

\begin{corollary}\label{cor-RHom-ABCD}
For $A,B\in\Dbc_\cx(X)$ and $C,D\in\Dbc_\cy(Y)$ there is
a functorial isomorphism
\[ \RHom(A,B)\tens_\CC \RHom(C,D) \rTTo_\sim
\RHom(A\dlgn C,B\dlgn D). \]
\end{corollary}

\begin{proof}
Follows from \lemref{lem-RGamma} with $E=\HOM(A,B)$
and $F=\HOM(C,D)$ and \corref{cor-RHom-box-isom}.
\end{proof}

Applying the K\"unneth formula we get

\begin{corollary}\label{cor-Ext2-Ext}
For $\cx$\coco\ $A,B$ and $\cy$\coco\ $C,D$ we have an isomorphism
\[ \oplus_{i+j=k} \Ext^i_{D(X)}(A,B) \tens_\CC
\Ext^j_{D(Y)}(C,D) \rTTo_\sim
\Ext^k_{D(X\times Y)}(A\dlgn C,B\dlgn D). \]
\end{corollary}

\begin{corollary}\label{cor-Hom2-Hom}
If $A,B\in\Perv(X)$, $C,D\in\Perv(Y)$, then
\begin{equation}
\Hom_{\Perv(X)}(A,B) \tens_\CC \Hom_{\Perv(Y)}(C,D)
\underset\sim\to \Hom_{\Perv(X\times Y)}(A\dlgn C,B\dlgn D).
\label{eq-cor-Hom}
\end{equation}
\end{corollary}

Indeed, for perverse sheaves $\Ext_{D(X)}^i(A,B)=0$ for $i<0$ by
\cite[Corollaire 2.1.4]{BeilBD:Perverse}.

\section{Deligne's external tensor product of perverse sheaves}
\label{sec-thm-main}
It follows from \cite{BeilBD:Perverse}, Amplification 1.4.17.1
and Proposition 1.4.18, by induction on strata that any
object of the abelian category $\Perv(X)$ has finite length.
Constructibility of perverse sheaves implies that the
$\CC$\n-vector spaces $\Hom_{\Perv(X)}(A,B)$ are finite
dimensional. Therefore, $\Perv(X)$ is an inductive limit of
its full subcategories $\<M\>$, equivalent to $A\modul$ for
some finite dimensional associative unital algebra $A$
\cite[Corollaire 2.17]{Del:cat}. The subcategory $\<M\>$
is formed by subquotients of $M^n$, $n\in\ZZ_{>0}$.
By results of Deligne
\cite{Del:cat} there exists an (abstract) external tensor
product of the categories $\Perv(X)$ and $\Perv(Y)$ --
the universal functor
\[ \Ddlgn: \Perv(X) \times \Perv(Y) \to
\Perv(X) \Ddlgn \Perv(Y), \]
whose target is some $\CC$\n-linear abelian category.
Universality implies, in particular, that there exists an
exact $\CC$\n-linear functor $F$ and an isomorphism
\begin{diagram}[LaTeXeqno]
\Perv(X) \times \Perv(Y) & \rTTo^\Ddlgn &
\Perv(X) \Ddlgn \Perv(Y) \\
& \rdTTo_\dlgn & \dTTo>F \\
&& \Perv(X\times Y)
\label{dia-Dlgn-T}
\end{diagram}
Our goal is to prove that $F$ is an equivalence. Once this is
done, we can choose the Deligne external tensor product
$\Ddlgn$ to be the geometric external tensor product
$\dlgn$ and $\Perv(X) \Ddlgn \Perv(Y)$ to be $\Perv(X\times Y)$.
Thus, we can use the same notation
$\dlgn$ in both abstract and geometric senses.

\begin{theorem}
The functor
\[ F: \Perv(X) \Ddlgn \Perv(Y) \to \Perv(X\times Y) \]
determined by diagram \eqref{dia-Dlgn-T} is an equivalence.
Therefore, we can choose as tensor product
$\Perv(X) \Ddlgn \Perv(Y) = \Perv(X\times Y)$ and
$\Ddlgn = \dlgn$.
\end{theorem}

\begin{proof}
The results of \secref{sec-simple-pvsh} describe
$\SimpPerv(X)$ -- the list of isomorphism classes of simple
objects of $\Perv(X)$. For $\Perv(X)\Ddlgn\Perv(Y)$ this
list is the tensor product $\SimpPerv(X)\Ddlgn\SimpPerv(Y)$
of the lists for $X$ and $Y$
\cite[Lemme 5.9]{Del:cat} due to algebraic
closedness of $\CC$. On the other hand, \propref{pro-6boxes}
implies that
\[ \SimpPerv(X\times Y) = \SimpPerv(X)\dlgn \SimpPerv(Y). \]
Therefore, the functor $F$ maps bijectively the list of
isomorphism classes of simple objects of
$\Perv(X)\Ddlgn\Perv(Y)$ to the list of isomorphism classes
of simple perverse sheaves of $\Perv(X\times Y)$.

\begin{lemma}[\cite{BeiGinSoe-Kosz_pat} Lemma 3.2.4]\label{lem-Ext-Hom}
The natural mapping of Yoneda's $\Ext$ to the $\Ext$
in the derived category
\[ \theta: \ExtY^i_{\Perv(X)}(A,B) \to \Ext^i_{D(X)}(A,B)
\overset{\text{def}}= \Hom_{D(X)}(A,B[i]) \]
is bijective for $k=0,1$ and injective for $k=2$
for all $A,B\in\Perv(X)$.
\end{lemma}

\begin{lemma}\label{lem-ExtY2-ExtY}
Let $\ca$, $\cb$ be $\CC$\n-linear abelian categories with
length (that is objects have finite length
and $\Hom$ spaces are finite dimensional). Then
\[ \oplus_{i+j=k} \ExtY^i_\ca(K,L) \tens_\CC \ExtY^j_\cb(M,N)
\rTTo^\sim \ExtY^k_{\ca\Ddlgn\cb}(K\Ddlgn M,L\Ddlgn N). \]
\end{lemma}

\begin{proof}
From the definition of Yoneda's $\ExtY$ \cite{Yoneda-Ext}
it is clear that
\begin{equation}
\ExtY^i_{\Perv(X)}(A,B) =
\lim_{\substack{\longrightarrow\\ \Perv(X)\supset\<P\>\ni A,B}}
\ExtY^i_{\<P\>}(A,B),
\label{eq-lim-Ext}
\end{equation}
where $P$ runs over such objects of $\Perv(X)$ that the
subcategory $\<P\>$ contains $A$ and $B$. Since the
category $\<P\>$ has enough injectives and projectives, we
can identify $\ExtY^i_{\<P\>}(C,D)$ with the right derived
functor $\Ext^i_{\<P\>}(C,D)$ of $\Hom_{\<P\>}$.

It suffices to prove the statement for subcategories of $\ca$,
$\cb$ of the form $\<P\>$. So we have to show that the
external product map
\[ \oplus_{i+j=k} \Ext^i_A(K,L) \tens_\CC \Ext^j_B(M,N)
\rTTo \Ext^k_{A\tens_\CC B}(K\tens_\CC M,L\tens_\CC N) \]
is bijective for finite dimensional $\CC$\n-algebras $A,B$,
finite dimensional $A$\n-modules $K$, $L$ and finite
dimensional $B$\n-modules $M$, $N$. This is precisely one of
the assertions of Theorem XI.3.1 of Cartan and Eilenberg
\cite{CarEil-book}.
\end{proof}

Combining the above lemmas and \corref{cor-Ext2-Ext} we
get a commutative diagram
\begin{diagram}[height=2.4em,nobalance,objectstyle=\scriptstyle]
\ExtY^k_{\Perv(X)\Ddlgn\Perv(Y)}(A\Ddlgn C,B\Ddlgn D) && \\
\dTTo<{\text{\lemref{lem-ExtY2-ExtY}}}>\wr && \\
\underset{i+j=k}{\bigoplus} \ExtY^i_{\Perv(X)}(A,B) \tens_\CC
\ExtY^j_{\Perv(Y)}(C,D) & \rTTo^{\theta} &
\underset{i+j=k}{\bigoplus} \Ext^i_{D(X)}(A,B) \tens_\CC \Ext^j_{D(Y)}(C,D) \\
\dTTo<{\mathrm V} &&
\dTTo<\wr>{\text{\corref{cor-Ext2-Ext}}} \\
\ExtY^k_{\Perv(X\times Y)}(A\dlgn C,B\dlgn D) & \rTTo^{\theta} &
\Ext^k_{D(X\times Y)}(A\dlgn C,B\dlgn D)
\end{diagram}
for all $A,B\in\Perv(X)$, $C,D\in\Perv(Y)$ and all
$k\in\ZZ_{\ge0}$. Here the map $\mathrm V$ takes
\begin{align*}
&[0\shortto B\shortto M_1\shortto\dots\shortto M_i\shortto A\shortto0] \tens
[0\shortto D\shortto N_1\shortto\dots\shortto N_j\shortto C\shortto0] \\
\intertext{to}
&[0\to B\dlgn D\to M_1\dlgn D\to\dots\to M_i\dlgn D\to \\
&\hspace*{47mm} \to A\dlgn N_1\to\dots\to A\dlgn N_j\to A\dlgn C\to0]
\end{align*}
where the middle map is the composition of
$M_i\dlgn D\to A\dlgn D \to A\dlgn N_1$. Commutativity of the
square follows from our sign convention. It is verified
similarly to Yoneda's computation
\cite{Yoneda-Ext-V-prod} of the $\mathsf V$ multiplication
of Cartan and Eilenberg \cite[Section XI.1]{CarEil-book}.

By \lemref{lem-Ext-Hom} the horizontal arrows are
bijective for $k=0,1$ and injective for $k=2$. Hence, the same
holds for the mapping ${\mathrm V}$. In particular, the
induced by $F$ mapping of Yoneda's $\ExtY^k$ between simple
objects of $\Perv(X)\Ddlgn\Perv(Y)$ to that of
$\Perv(X\times Y)$ is bijective for $k=0,1$ and injective for
$k=2$.  It remains to apply the following

\begin{lemma}\label{lem-Ext-01bij-2inj}
Let $F:\ca\to\cb$ be an exact functor between essentially
small categories with length. Assume that $F$
induces a bijection on the list of isomorphism classes of
simple objects. Assume also that the maps induced by $F$
\begin{equation}
\ExtY^k_\ca(T,S) \to \ExtY^k_\cb(FT,FS)
\label{eq-ExtY-TS}
\end{equation}
are bijective for $k=0,1$ and injective for $k=2$ for all
simple objects $T$, $S$ of $\ca$. Then $F$ is an equivalence.
\end{lemma}

\begin{proof}
First we prove by induction on the length of $T$ that
\eqref{eq-ExtY-TS} is an isomorphism for $k=0,1$ for a
simple $S$ and an arbitrary $T$. Indeed, write the long
exact sequences up to $k=2$ for $\ExtY^\bull_\ca(\text-,S)$
and $\ExtY^\bull_\cb(F\text-,FS)$ associated with
$0\to T'\to T\to T''\to0$, where $T''$ is simple,  and use the
5\n-Lemma.
Second, we prove by induction on the length of $S$ that
\eqref{eq-ExtY-TS} is an isomorphism for $k=0$ for all
objects $S$, $T$ of $\ca$. Indeed, write the long
exact sequences up to $k=1$ for $\ExtY^\bull_\ca(T,\text-)$
and $\ExtY^\bull_\cb(FT,F\text-)$ associated with
$0\to S'\to S\to S''\to0$, where $S'$ is simple,  and use the
5\n-Lemma. Hence, $F$ is full and faithful.

Now we prove by induction on length that $F$ induces
surjection on the set of isomorphism classes of objects of
length $\le n$ for $\ca$ and $\cb$. For $n=1$ it is the
hypothesis of the lemma. For $X\in\Ob\cb$ of length $n>1$
we can assume the existence of
\begin{equation}
0 \to FS \to X \to FT \to 0
\label{eq-FS-X-FT}
\end{equation}
for some simple $S\in\Ob\ca$ and some $T\in\Ob\ca$ of
length less than $n$. Since the map
\[ \ExtY^1_\ca(T,S) \to \ExtY^1_\cb(FT,FS) \]
is bijective, there exists a short exact sequence
$0 \to S \to Y \to T \to 0$ in $\ca$ such that
$0 \to FS \to FY \to FT \to 0$ is congruent with
\eqref{eq-FS-X-FT}. In particular, $X\simeq FY$.

Therefore, $F$ is full, faithful and essentially surjective on objects.
\end{proof}

Applying this lemma to $\ca=\Perv(X)\Ddlgn\Perv(Y)$ and
$\cb=\Perv(X\times Y)$ we prove the theorem.
\end{proof}

\appendix
\section{Some formulas}

Here we summarise some formulas taken mostly
from Borel \cite{Borel-Int_Cohom-V} \S10. All spaces are
locally compact, locally completely paracompact,
locally contractible and of finite cohomological
dimension over $\CC$.

For a continuous map $f:X\to Y$ and $A, B\in D(Y)$ we have
\begin{equation}
f^*(A\tens B) \simeq f^*A\tens f^*B,
\label{ApP10.1}
\end{equation}
by \lococit\ Proposition 10.1.

For $A, B, C\in D(X)$ we have
\begin{equation}
\RHOM(A\tens B,C) \simeq \RHOM(A,\RHOM(B,C))
\label{ApP10.2}
\end{equation}
by \lococit\ Proposition 10.2.

For a continuous map $f:X\to Y$ and $A\in D(Y)$,
$B\in D(X)$ we have
\begin{equation}
f_*\RHOM(f^*A,B) \simeq \RHOM(A,f_*B),
\label{ApP10.3(1)}
\end{equation}
by \lococit\ Proposition 10.3(1) and
\begin{equation}
f_!(B\tens f^*A) \simeq f_!B\tens A,
\label{ApP10.8(2)}
\end{equation}
by \lococit\ Proposition 10.8(2).

Denoting $p_1:X\times Y\to X$ the projection
on the first space, we have for $A,B\in D^b(X)$
\begin{equation}
p_1^*\RHOM(A,B) \simeq \RHOM(p_1^*A,p_1^*B)
\label{ApP10.21}
\end{equation}
by \lococit\ Proposition 10.21.

Denoting $p_2:X\times Y\to Y$ the projection
on the second space, we have for $\cx$\clc
$A\in D^b_\cx(X)$ and $\cy$\coco\ $B\in\Dbc_\cy(Y)$
\begin{equation}
p_1^*\D A\tens p_2^*B \simeq \RHOM(p_1^*A,p_2^!B)
\label{ApT10.25}
\end{equation}
by \lococit\ Theorem 10.25.

Here
\[ \D A = \RHOM(A,\D_X) \in D^b_\cx(X) \]
is the Verdier dual of $A$ and
\[ \D_X = g^!\CC \in\Dbc_\cx(X), \qquad g:X\to pt, \]
is the dualising sheaf \cite{Verdier1} (see also
\cite[7.18 and Theorem 8.3]{Borel-Int_Cohom-V}).
Substituting $A=\CC$ into \eqref{ApT10.25} one gets
\begin{equation}
\D_X\dlgn B \simeq p_2^!B.
\label{eq-Ap-p!B}
\end{equation}

We have
\begin{equation}
\D_X\dlgn\D_Y\simeq\D_{X\times Y}
\label{ApC10.26}
\end{equation}
by \cite{Borel-Int_Cohom-V} Corollary 10.26, and
\[ \D^2 \simeq \Id \]
by \cite{Verdier1} (see also \cite{Borel-Int_Cohom-V} Theorem 8.10).
If $f:X\to Y$ is a stratified map, then
\begin{equation}
\D f^* = f^!\D, \qquad \D f^! = f^*\D,
\label{ApT10.17(1)}
\end{equation}
if the stratified map $f$ is proper, or algebraic over
$\CC$, or if every fibre of $f$ is compactifiable, then
\begin{equation}
\D f_* = f_!\D, \qquad \D f_! = f_*\D
\label{ApT10.17(2)}
\end{equation}
by \cite{Borel-Int_Cohom-V} Theorem 10.17.

\begin{definition}[see \cite{BeilBD:Perverse} Proposition 2.1.9]
\label{def-Ap-j!*}
Let $U\subset X$ be an open stratified subspace of
a stratified pseudomanifold $X=(X,\cx,p)$.
Denote by $j:U \rMono X$ the inclusion.
Let $A\in\Perv(U,U\cap\cx,p)$.
Its \emph{prolongation} $j_{!*}A$ is an object
$B\in\Dbc_\cx(X)$ equipped with an isomorphism
$j^*B\simeq A$, such that for any stratum
$s:S \rMono X$, that is not contained in $U$,
we have $H^is^*B=0$ for $i\ge p(S)$ and $H^is^!B=0$ for $i\le p(S)$.
The object $B$ is determined uniquely up to a
unique isomorphism. This gives a functor
$j_{!*}:\Perv(U,U\cap\cx,p)\to\Perv(X,\cx,p)$.
\end{definition}

If $k:V \rMono U$ is another stratified open inclusion,
then there is an isomorphism (2.1.7.1) of \cite{BeilBD:Perverse}
\begin{equation}
(j\circ k)_{!*} \rTTo^\sim j_{!*}\circ k_{!*}.
\label{Ap2.1.7.1}
\end{equation}

%\bibliography{yuri}
%\end{document}

\end{document}